# Uniform deterministic equivalent of additive functionals and non-parametric drift estimation for one-dimensional recurrent diffusions

D. Loukianova[a] and O. Loukianov[b]

[a]*Département Mathématique, Université d'Evry, France.* E-mail: dasha.loukianova@univ-evry.fr
[b]*Département Informatique, IUT de Fontainebleau, France.* E-mail: oleg@iut-fbleau.fr

**Abstract.** Usually the problem of drift estimation for a diffusion process is considered under the hypothesis of ergodicity. It is less often considered under the hypothesis of null-recurrence, simply because there are fewer limit theorems and existing ones do not apply to the whole null-recurrent class.

The aim of this paper is to provide some limit theorems for additive functionals and martingales of a general (ergodic or null) recurrent diffusion which would allow us to have a somewhat unified approach to the problem of non-parametric kernel drift estimation in the one-dimensional recurrent case. As a particular example we obtain the rate of convergence of the Nadaraya–Watson estimator in the case of a locally Hölder-continuous drift.

**Résumé.** Habituellement le problème de l'estimation du drift pour un processus de diffusion est considéré sous l'hypothèse de l'ergodicité. Il l'est moins souvent sous l'hypothèse de nulle-récurrence, car dans ce cas il y a moins de théorèmes limites, et ceux qui existent ne s'appliquent pas à toute la classe nulle-récurrente.

Le but de cet article est de démontrer quelques théorèmes limites pour les fonctionnelles additives et martingales dépendantes d'une diffusion récurrente générale (ergodique ou nulle). Ces théorèmes permettent de donner une approche unifiée au problème de l'estimation non-paramétrique par noyau du drift dans le cas unidimensionnel récurrent. Comme exemple on obtient la vitesse de convergence de l'estimateur de Nadaraya–Watson dans le cas d'un drift localement hölderien.

## 1. Introduction

Consider a stochastic differential equation

$$dX_t = \sigma(X_t)\,dW_t + b(X_t)\,dt, \tag{1}$$

where $b, \sigma: \mathbb{R} \to \mathbb{R}$ and $(W_t)_{t \geq 0}$ is a Brownian motion on $\mathbb{R}$. Throughout this paper we suppose that $\sigma$ is strictly positive, continuous, $b$ is measurable, and for some constant $C > 0$ and all $x \in \mathbb{R}$

$$\sigma^2(x) \leq C(1 + x^2) \quad \text{and} \quad |b(x)| \leq C(1 + |x|).$$

Under these conditions for each $x_0 \in \mathbb{R}$ the Eq. (1) has a unique weak solution $X = (X_t)_{t \geq 0}$ starting from $x_0$, and the corresponding semigroup is strong Feller (see e.g. [20], p. 170). Furthermore, we suppose that $X$ is Harris recurrent. Recall that it means that $X$ admits an invariant measure $\mu$ on $\mathcal{B}(\mathbb{R})$ such that for any measurable $f \geq 0$,

$$\mu(f) > 0 \implies \forall x \in \mathbb{R}: \int_0^\infty f(X_t) \, dt = \infty \quad P_x\text{-a.s.}$$

In our case (more exactly if $\sigma$ and $b$ are locally bounded Borel functions and $\sigma$ does not vanish), the invariant measure $\mu$ is absolutely continuous with respect to the Lebesgue measure, see [19], p. 298. The density of $\mu(dx)$ is given by

$$\mu(dx) = \frac{2}{\sigma^2(x)} \exp\left(\int_0^x \frac{2b}{\sigma^2}(v) \, dv\right) dx. \tag{2}$$

Recall also that a one-dimensional diffusion given by (1) is recurrent if and only if

$$\int_0^x \exp\left(\int_0^y \frac{2b}{\sigma^2}(v) \, dv\right) dy \to \pm\infty \quad \text{as } x \to \pm\infty. \tag{3}$$

With regard to the $\mu$-total mass of $\mathbb{R}$ one subdivides Harris recurrent diffusions into two sub-classes: if $\mu(\mathbb{R})$ is finite then $X$ is said to be ergodic, and null-recurrent otherwise.

For such a diffusion we will consider the problem of drift estimation. Throughout this paper we suppose that a sample path of diffusion is observed on $[0, T]$ and that $T \to \infty$.

Since we study the asymptotic behavior of diffusion processes, it is important to have at one's disposal limit theorems for additive functionals of the form $A_t = \int_0^t f(X_s) \, ds$ or martingales of the form $M_t = \int_0^t f(X_s) \, dW_s$. The limit theorems for additive functionals and martingales of *ergodic* diffusions are available on each level of convergence and well known, for instance if $f \geq 0$, $\mu(f) < \infty$, we have a LLN

$$\lim_{t \to \infty} \frac{1}{t} \int_0^t f(X_s) \, ds = \mu(f) \quad P_x\text{-a.s. } \forall x$$

and if $\mu(f^2) < \infty$, a CLT:

$$\frac{1}{\sqrt{t}} \int_0^t f(X_s) \, dW_s \longrightarrow \mathcal{N}(0, \mu(f^2)).$$

For *null-recurrent* diffusions there is no LLN with deterministic normalization and weak convergence theorems are available only for a subclass of regularly varying null-recurrent diffusions, see [12, 21] and the book by Höpfner and Löcherbach [10]. As for limit theorems working in a general recurrent setting, there is only one result: the famous Chacon–Ornstein theorem stating that all integrable additive functionals (IAF) of $X$ are asymptotically equivalent: $\forall f \geq 0$, $\mu(f) < \infty$; $\forall g \geq 0$, $0 < \mu(g) < \infty$;

$$\lim_{t \to \infty} \frac{\int_0^t f(X_s) \, ds}{\int_0^t g(X_s) \, ds} = \frac{\mu(f)}{\mu(g)} \quad P_x\text{-a.s. } \forall x$$

and its integral version yields that $\mu$-a.s.

$$\lim_{t \to \infty} \frac{\mathbb{E}_x \int_0^t f(X_s) \, ds}{\mathbb{E}_x \int_0^t g(X_s) \, ds} = \frac{\mu(f)}{\mu(g)}; \tag{4}$$

where the exceptional set depends on $f$ and on $g$.

The literature on statistical inference for ergodic diffusions in general and drift estimation in particular is extensive and significant results can be obtained in this case, see for example the works by Dalalyan and Kutoyants [5], Dalalyan [4], Van-Zanten [22], Galtchouk and Pergamentchikov [8], Yoshida [24] to mention just a few; an extensive survey can be found in the recent book [14] by Kutoyants.

The situation is different for null-recurrent diffusions. As we have mentioned, an important class of *regularly varying* null-recurrent diffusions behaving much like the ergodic ones has been thoroughly studied by Höpfner and Löcherbach in [10]. For such diffusions the pair $(M_t/\sqrt{v_t}, A_t/v_t)$ converges in law as $t \to \infty$, where $v_t = t^\alpha l(t)$ for some $0 < \alpha \leq 1$ and $l(t)$ varying slowly at infinity. Using these facts, Höpfner and Kutoyants [9] have presented one of the first examples concerning the rate of convergence and the limit distribution of MLE for null-recurrent diffusions. It might be possible to extend their method to the whole class of regular variation, but since it is based on weak convergence it cannot be used in general cases.

Our aim in this paper is to develop some tools, for the problem of drift estimation, and more precisely, for the rate of convergence calculation, which would work in a general recurrent (null or ergodic) setting.

The first idea in this direction could be derived from the Chacon–Ornstein theorem. Because all IAF of $X$ are equivalent, we can consider some fixed IAF as a "time" of the system and try to work with a random normalisation based on this IAF. This idea was used in [7] to study the rate of convergence of a non-parametric kernel estimator and by Loukianova, Loukianov [17] to obtain the a.s. rate of convergence of the MLE. However, applying this idea in each particular case presents a real technical difficulty, and no concise method, based on this idea, has emerged so far.

The second idea is based on the following observation: we are interested in the calculation of the rate of convergence, and the rate is a notion of tightness. All known limit theorems deal with a stronger type of convergence, and none (except Chacon–Ornstein) works for the whole recurrent class. However, tightness should be sufficient to treat the problem of the rate of convergence. Therefore the natural question arises: is there a property of tightness with deterministic norming for IAF, which would be true for all recurrent diffusions? The answer is affirmative:

For every recurrent diffusion $X$ there is some deterministic function $v_t$, called in the sequel *deterministic equivalent of $X$*. With this function, which will be explained shortly, the following theorem holds:

**Theorem 1.** *For every IAF $A_t$ of $X$ the following holds*

$$\lim_{M \to \infty} \liminf_{t \to \infty} P\left(\frac{1}{M} < \frac{A_t}{v_t} < M\right) = 1.$$

Such a property was first proved for recurrent Markov chains by Chen [3], without any relation to statistics. In [16], using Chen's method, we extended this property to one-dimensional diffusions, and recently, in [15], using a new version of the Nummelin splitting method, this property was extended to the whole class of continuous time Harris recurrent Markov processes.

To explain what a deterministic equivalent $v_t$ is, we need the notion of a special function.

**Definition 1.1.** *A measurable bounded function $f : \mathbb{R} \to \mathbb{R}^+$ is said to be special if for every $h : \mathbb{R} \to \mathbb{R}^+$, measurable, bounded and such that $\mu(h) > 0$, the following function $U^h f$ is bounded:*

$$U^h f(x) = \mathbb{E}_x \int_0^\infty \exp\left[-\int_0^t h(X_s)\,ds\right] f(X_t)\,dt.$$

This notion can be extended to additive functionals.

**Definition 1.2.** *A continuous AF $A_t$ is said to be special if for every $h : \mathbb{R} \to \mathbb{R}^+$, measurable, bounded and such that $\mu(h) > 0$, the function*

$$U_A^h 1(x) = \mathbb{E}_x \int_0^\infty \exp\left[-\int_0^t h(X_s)\,ds\right] dA_t$$

*is bounded.*

For instance, a diffusion local time is a special additive functional (SAF), see [2]; on the other hand, any measurable bounded function with compact support is special for a $n$-dimensional strong Feller diffusion, see [15].

Observe that in the "integral version" (4) of the Chacon–Ornstein theorem we cannot determine whether the result holds or not for fixed $x$, $f$ and $g$. The most important application of special functions and functionals is in the fact that the assertion (4) holds for every $x$, and moreover, the following "strong" version of Chacon–Ornstein theorem (SCO) holds:

**Theorem 2.** *If $A$, $B$ are two SAF such that $\|\nu_B\| := \mathbb{E}_\mu B_1 > 0$, then for every pair $(\pi_1, \pi_2)$ of probability measures,*

$$\lim_{t \to \infty} \frac{\mathbb{E}_{\pi_1} A_t}{\mathbb{E}_{\pi_2} B_t} = \frac{\|\nu_A\|}{\|\nu_B\|}.$$

(Recall that if $A_t = \int_0^t f(X_s)\,ds$ with $f \in L^1(\mu)$, then $\|\nu_A\| = \mu(f)$.)

Now, to normalize our additive functional we take some special $g$ with $\mu(g) = 1$, and put

$$v_t = \mathbb{E}_\pi \int_0^t g(X_s)\,ds,$$

where $\pi$ is some probability measure on $E$. Clearly, $v_t$ is non-negative and non-decreasing. In view of the strong Chacon–Ornstein theorem the asymptotic order of $v_t$ depends only on the law of $X$. If $X$ is an ergodic diffusion, then $v_t \sim t$, if $X$ is a null-recurrent diffusion of regular variation with index $\alpha$, then $v_t \sim t^\alpha$; in particular, if $X$ is a linear Brownian motion, then $v_t \sim \sqrt{t}$. Some less explicit examples are possible: if $X$ is a drifted Brownian motion with compactly supported drift, then $X$ is null-recurrent and $v_t \sim \sqrt{t}$. However, to evaluate $v_t$ in the general case we need the semi-group of $X$, so $v_t$ is clearly a more "abstract" object.

In [16] we have applied the deterministic equivalent to give a rate of convergence of the MLE for a class of recurrent diffusions with Hölder drift. In this parametric context, the deterministic equivalent turned out to be an appropriate tool, because if in the proof of the ergodic case, we replace the normalization $t$ with $v_t$, we obtain, with very few modifications, a proof for the whole recurrent class.

In this paper we apply the deterministic equivalent to the calculation of the rate of convergence in non-parametric drift estimation. We suppose, as usual, that the drift function of $X$ is locally $\alpha$-Hölder near some point $x_0 \in E$, and we estimate $b(x_0)$ by the Nadaraya–Watson estimator $\hat{b}^{h_t}_{x_0,t}$ given by

$$\hat{b}^h_{x_0,t} = \frac{\int_0^t \phi((X_s - x_0)/h)\,dX_s}{\int_0^t \phi((X_s - x_0)/h)\,ds}.$$

Here $\phi$ is a non-negative smooth function with compact support. The process $\hat{b}^h_{x_0,t}$ is composed with some bandwidth $h_t$. In the ergodic setting $h_t$ is a deterministic function of $t$ which can be found using a standard optimization procedure. A similar procedure, when guessing the optimal bandwidth in the general recurrent case, gives $h_t$ as a function of $v_t$. As we have noticed above, $v_t$ cannot be calculated in general, but we need an estimator, depending only on the observations. So we also express the bandwidth (and the rate) in terms of any random and observable IAF $V_t$, equivalent (in the sense of Theorem 1) to $v_t$.

Therefore, the objects appearing in this problem and those for which we need to study the asymptotic behaviour, are additive functionals $A^h_t$ and martingales $M^h_t$ composed with some deterministic or random bandwidth $h_t$, namely

$$A^{h_t}_t = \int_0^t \phi\left(\frac{X_s - x_0}{h_t}\right) ds \quad \text{and} \quad M^{h_t}_t = \int_0^t \phi\left(\frac{X_s - x_0}{h_t}\right) dW_s.$$

To be precise, with a view of rate calculation, we want to show some tightness property for these objects, with some normalization depending on $v_t$. When the bandwidth is random, $A_t^{h_t}$ is no longer an additive functional, and $M_t^{h_t}$ is no longer a martingale. It would be difficult to obtain their asymptotics directly, using stochastic calculus methods, because these objects do not have the usual stochastic structure. To avoid this difficulty, we firstly show some uniform in $h$ tightness property for $A_t^h$. Namely, we show that $\frac{1}{h}A_t^h$ admits $v_t$ as uniform deterministic equivalent (Theorem 5). We then derive a "diagonal" tightness property for $A_t^{h_t}$ and $M_t^{h_t}$, Corollaries 3.1 and 3.2. All results regarding the limit behavior of $A_t^{h_t}$ and $M_t^{h_t}$ are given in Section 3.

The main technical tools of our work are two theorems on the local time $L_t^x$ of $X$, namely, the uniform in $x$ deterministic equivalent (Theorem 3), and the uniform in $x$ strong Chacon–Ornstein Theorem 4. This is the content of Section 2.

Statistical applications are given in Section 4. Theorem 6 gives the rate $r_t = v_t^{-\alpha/(2\alpha+1)}$ specified by the bandwidth $h_t = v_t^{-1/(2\alpha+1)}$. The random rate $R_t$ and the bandwidth $H_t$ are given by the same theorem and are: $R_t = V_t^{-\alpha/(2\alpha+1)}$; $H_t = V_t^{-1/(2\alpha+1)}$.

Expressed in deterministic terms, our result agrees with the well-known rate $r_t = t^{-\alpha/(2\alpha+1)}$ in the ergodic case of this model, as well as with those of Delattre and Hoffmann [6], where the minimax rate $r_t = t^{-\alpha/(4\alpha+2)}$ specified by the bandwidth $h_t = t^{-(-1/(4\alpha+2))}$ was found for the model of drifted Brownian motion with compactly supported drift. Recall that this diffusion is null-recurrent with $v_t \sim \sqrt{t}$. Our random expressions agree with those of Delattre, Hoffmann and Kessler [7], where the Nadaraya–Watson estimator was studied in a minimax setting by different methods and both rate and bandwidth were expressed in terms of one fixed IAF: the occupation time.

Notice finally that for multi-dimensional diffusions, it is possible, using the Nummelin splitting, to get a result similar to Theorem 6, but only with deterministic rate and bandwidth [15].

## 2. Uniform deterministic equivalent and uniform Strong Chacon–Ornstein theorem for the local time

Recall that we consider a recurrent scalar diffusion $X$ given by the Eq. (1) and subject to the same assumptions. As usual, $\mu$ denotes the invariant measure of $X$ and $(L_t^x)_{t \geq 0}$ the Tanaka–Meyer local time at the point $x \in \mathbb{R}$. Let $\pi$ be some probability on $\mathbb{R}$ and $v_t = \mathbb{E}_\pi \int_0^t g(X_s)\,ds$ be a deterministic equivalent of $X$. We suppose $g$ measurable, bounded, compactly supported and such that $\mu(g) = 1$. In this section we derive a uniform deterministic equivalent for $(L_t^x)_{t \geq 0}$ (Theorem 3), and we show that the local time obeys a uniform strong Chacon–Ornstein theorem (Theorem 4).

**Proposition 2.1.** *For each $y \in \mathbb{R}$, $\|\nu_{L^y}\| = \sigma^2(y)\mu(y)$, where $\mu(y)$ is the density of the invariant measure $\mu(dy)$ with respect to Lebesgue measure. In particular, $0 < \|\nu_{L^y}\| < \infty$.*

**Proof.** Denote by $p(t; x, y)$ the density of the transition function of $X$ with respect to the invariant measure $\mu$. This density exists and may be taken to be positive, jointly continuous in all variables and symmetric, that is $p(t, x, y) = p(t, y, x)$, see [11], p. 149.

On the other hand, the invariant measure $\mu$ has itself a density with respect to Lebesgue measure, $\mu(dx) = \mu(x)\,dx$, where

$$\mu(x) = \frac{2}{\sigma^2(x)} \exp\left(\int_0^x \frac{2b}{\sigma^2}(v)\,dv\right). \tag{5}$$

We have the following relation between the local time and the density $p(t; x, y)$:

$$\mathbb{E}_x L_t^y = \sigma^2(y)\mu(y) \int_0^t p(s; x, y)\,ds, \tag{6}$$

see [1], p. 21 or [11], p. 175. The factor $\sigma^2(y)\mu(y)$ arises since $L_t^y$ is the Tanaka–Meyer local time and not the Ito–McKean local time.

$$\|\nu_{L^y}\| := \mathbb{E}_\mu L_1^y = \sigma^2(y)\mu(y) \int_\mathbb{R} \mu(x)\,dx \int_0^1 p(s;x,y)\,ds$$
$$= \sigma^2(y)\mu(y) \int_0^1 ds \int_\mathbb{R} p(s;y,x)\mu(x)\,dx = \sigma^2(y)\mu(y).$$

$\square$

**Remark 2.1.** *The condition $\|\nu_{L^y}\| > 0$ implies, according to Brancovan [2], that $\forall x \in \mathbb{R}$, $\mathbf{P}_x(L_\infty^y = \infty) = 1$, and that $L^y$ is a SAF.*

**Theorem 3.** *Let $K \subset \mathbb{R}$ be some compact set. For any initial probability $\nu$ the following limit holds:*

$$\lim_{M\to\infty} \liminf_{t\to\infty} \mathbf{P}_\nu\left(\frac{1}{M} \leq \inf_{x\in K} \frac{L_t^x}{v_t} \leq \sup_{x\in K} \frac{L_t^x}{v_t} \leq M\right) = 1.$$

Before we prove this theorem let us discuss the structure of the proof. The basic idea is of course to show some kind of uniform continuity property of the family of normalized local times. For each point $x_0 \in K$, we do it for some small interval $[x_0 - \delta; x_0 + \delta]$ with $\delta > 0$ depending on $x_0$. Using this we will start by proving the Theorem 3 for this small interval and then for the compact $K$. By definition for $y < z$,

$$\frac{L_t^y - L_t^z}{2} = (X_t - y)^+ - (X_t - z)^+ - [(X_0 - y)^+ - (X_0 - z)^+]$$
$$- \int_0^t \mathbf{1}_{\{y<X_s\leq z\}}\sigma(X_s)\,dW_s - \int_0^t \mathbf{1}_{\{y<X_s\leq z\}}b(X_s)\,ds \quad \text{a.s.} \tag{7}$$

For $t \geq 0$, $y < z$ denote

$$M_t^{y,z} = \int_0^t \mathbf{1}_{\{y<X_s\leq z\}}\sigma(X_s)\,dW_s. \tag{8}$$

Using occupation formula and boundedness of $|b|/\sigma^2$ on $K$, with $C > 0$, we can write

$$|L_t^y - L_t^z| \leq 2\delta + |M_t^{y,z}| + C\int_y^z L_t^x\,dx \quad \text{a.s.} \tag{9}$$

Now to prove Theorem 3 we need Lemmas 2.1–2.4. The first one is axillary to the second. The second deals with the martingale term of the decomposition (9), the third one with the last term of the (9) and Lemma 2.4 gives the equicontinuity property we need.

Put $K_\delta := [x_0 - \delta; x_0 + \delta]$.

**Lemma 2.1.** *For any initial probability $\nu$ the following limit holds:*

$$\lim_{M\to\infty} \limsup_{t\to\infty} \mathbf{P}_\nu\left(\int_{x_0-\delta}^{x_0+\delta}(L_t^x)^2\,dx > Mv_t^2\right) = 0.$$

**Proof.** Fix some $x < x_0$. Observe that with some constant $C > 0$ we have

$$(L_t^x - L_t^{x_0})^2 \leq C\left\{(2\delta)^2 + \left(\int_0^t \mathbf{1}_{\{x<X_s\leq x_0\}}\sigma(X_s)\,dW_s\right)^2 + \left(\int_0^t \mathbf{1}_{\{X_s\in K_\delta\}}|b(X_s)|\,ds\right)^2\right\}.$$

Observe that thanks to conditions on $b$ and $\sigma$ in the Introduction they are bounded on $K_\delta$. For $t \geq 0$, introduce the notation

$$A_t = \int_0^t \mathbf{1}_{\{X_s \in K_\delta\}} |b(X_s)| \, \mathrm{d}s. \tag{10}$$

Since $A_t$ is an IAF (and even SAF), with notation (10) we have

$$(L_t^x - L_t^{x_0})^2 \mathbf{1}_{A_t \leq \sqrt[4]{M} v_t} \leq C\{(2\delta)^2 + (M_t^{x,x_0})^2 + \sqrt{M} v_t^2\}.$$

Notice that

$$\mathbb{E}_\nu (M_t^{x,x_0})^2 \leq \mathbb{E}_\nu \int_0^t \mathbf{1}_{\{X_s \in K_\delta\}} \sigma^2(X_s) \, \mathrm{d}s := \tilde{v}_t.$$

Thanks to the condition on $\sigma$ the last expectation is an expectation of SAF, and we have denoted it by $\tilde{v}_t$ to stress the fact that it is a version of the deterministic equivalent of $X$. Observe that the same reasoning (with changing $M_t^{x,x_0}$ on $M_t^{x_0,x}$) and the same estimations hold for $x > x_0$. Finally, for all $x \in K_\delta$ it holds

$$\mathbb{E}_\nu (L_t^x - L_t^{x_0})^2 \mathbf{1}_{A_t \leq \sqrt[4]{M} v_t} \leq C\{(2\delta)^2 + \tilde{v}_t + \sqrt{M} v_t^2\}.$$

Now, to prove the lemma we use the following decomposition:

$$\mathbf{P}_\nu \left( \int_{x_0-\delta}^{x_0+\delta} (L_t^x)^2 \, \mathrm{d}x > M v_t^2 \right) \leq \mathbf{P}_\nu (8\delta (L_t^{x_0})^2 > M v_t^2) + \mathbf{P}_\nu \left( 4 \int_{x_0-\delta}^{x_0+\delta} (L_t^x - L_t^{x_0})^2 \, \mathrm{d}x > M v_t^2; A_t < \sqrt[4]{M} v_t \right)$$
$$+ \mathbf{P}_\nu (A_t > \sqrt[4]{M} v_t). \tag{11}$$

With some constant $C$, for the second term of the right-hand side expression of (11) we have:

$$\mathbf{P}_\nu \left( 4 \int_{x_0-\delta}^{x_0+\delta} (L_t^x - L_t^{x_0})^2 \, \mathrm{d}x > M v_t^2; A_t < \sqrt[4]{M} v_t \right)$$
$$\leq \frac{1}{M v_t^2} \int_{x_0-\delta}^{x_0+\delta} \mathbb{E}_\nu (L_t^x - L_t^{x_0})^2 \mathbf{1}_{A_t \leq \sqrt[4]{M} v_t} \, \mathrm{d}x \leq C \frac{(2\delta)^2 + \tilde{v}_t + \sqrt{M} v_t^2}{M v_t^2}.$$

According the strong Chacon–Ornstein theorem there is a limit of $\tilde{v}_t/v_t$, as $t \to \infty$. Using this, together with Theorem 1 for the other two terms on the right of (11) we obtain Lemma 2.1. □

Let $V$ be some IAF of $X$, such that $\|\nu_V\| > 0$, and for simplicity we take $\|\nu_V\| = 1$. $V_t$ will play the role of normalization for $L_t^y - L_t^z$ to obtain the uniform continuity property we need.

In the sequel of this section denote

$$\Omega_{M,t} := \left\{ \int_{K_\delta} (L_t^x)^2 \, \mathrm{d}x < M^2 v_t^2; \frac{v_t}{M} \leq V_t \leq v_t M \right\}.$$

Let $\{(M_t^{y,z}; t \geq 0); x_0 - \delta \leq y < z \leq x_0 + \delta\}$ be a family of martingales defined by (8).

**Lemma 2.2.** *For every $u > 0$ and $M > 0$ there is some $\tau = \tau(\omega, M, u, K_\delta)$ such that a.s. for $t > \tau$ it holds*

$$\sup_{y,z \in K_\delta} \frac{|M_t^{y,z}|}{V_t} \mathbf{1}_{\Omega_{M,t}} < \frac{u}{v_t^{1/4}}.$$

**Proof.** Fix some $u > 0$ and $M > 0$. Denote

$$B_n = \left\{ \sup_{\{t; e^{n-1} < v_t \leq e^n\}} \sup_{\{a,b \in K_\delta\}} \frac{|M_t^{y,z}|}{V_t} \mathbf{1}_{\Omega_{M,t}} \geq \frac{u}{v_t^{1/4}} \right\}.$$

We will use the Nishiyama martingale inequality to bound the probability of $B_n$ and then the Borel–Cantelli lemma to prove Lemma 2.2. The statement of the Nishiyama martingale's inequality for our setting is given bellow, for a proof of this result see [18] and for its generalization see [23].

Let $\rho$ be a metric on $K_\delta$, denote by $\|M\|_{\rho,t}$ the quadratic $\rho$-modulus of the family of martingales $\{M_t^{y,z}; y < z, y, z \in K_\delta\}$. For every $t > 0$, this modulus is given by

$$\|M\|_{\rho,t} = \sup_{y \neq z} \frac{\sqrt{\langle M_t^{y,z} \rangle}}{\rho(y,z)}.$$

As usual, denote by $N(K_\delta, \rho, \varepsilon)$ the metric entropy of $K_\delta$ with respect to the metric $\rho$, i.e. the minimal number of $\rho$-balls of radius $\varepsilon$, needed to cover $K_\delta$. The Nishiyama inequality for tail probabilities [23] states that

$$\mathbf{P}\left( \sup_{t \leq \tau} \sup_{\rho(y,z) \leq \eta} |M_t^{y,z}| \mathbf{1}_{\{\|M\|_{\rho,\tau} \leq \kappa\}} \geq x \right) \leq 2 e^{-x^2 / (c\phi^2(\eta)\kappa^2)},$$

where $C > 0, \eta > 0, \kappa > 0$ are constants, $\tau$ is a finite stopping time and

$$\phi(\eta) = \int_0^\eta \sqrt{\log N(K_\delta, \rho, \varepsilon)} \, d\varepsilon$$

is supposed to be finite.

Before we use the Nishiyama inequality we need to evaluate the quadratic $\rho$-modulus on the set $\Omega_{M,t}$. We have

$$\langle M_t^{y,z} \rangle = \int_0^t \mathbf{1}_{\{[y;z[}(X_s)\sigma^2(X_s) \, ds = \int_y^z L_t^x \, dx \leq \sqrt{z-y} \sqrt{\int_{x_0-\delta}^{x_0+\delta} (L_t^x)^2 \, dx}$$

and then

$$\langle M_t^{y,z} \rangle \mathbf{1}_{\Omega_{M,t}} \leq \sqrt{z-y} M v_t.$$

If we take $\rho(y,z) = \sqrt[4]{z-y}$ as a distance on $K_\delta$, it is easy to see that $\phi(\delta)$ is finite and we obtain a bound for quadratic $\rho$-modulus $\|M\|_{\rho,t}$ restricted on $\Omega_{M,t}$:

$$\|M\|_{\rho,t} \mathbf{1}_{\Omega_{M,t}} = \sup_{y \neq z} \frac{\sqrt{\langle M_t^{y,z} \rangle}}{\rho(y,z)} \mathbf{1}_{\Omega_{M,t}} \leq \sqrt{M v_t}.$$

Then

$$\mathbf{P}(B_n) = \mathbf{P}\left( \sup_{\{t; e^{n-1} < v_t \leq e^n\}} \sup_{\{y,z \in K_\delta\}} \frac{|M_t^{y,z}|}{V_t} \mathbf{1}_{\Omega_{M,t}} > \frac{u}{v_t^{1/4}} \right)$$

$$\leq \mathbf{P}\left( \sup_{\{t; e^{n-1} < v_t \leq e^n\}} \sup_{\{y,z \in K_\delta\}} |M_t^{y,z}| > \frac{u v_t}{M v_t^{1/4}}; \|M\|_{\rho,t} \leq \sqrt{M v_t} \right)$$

$$\leq \mathbf{P}\left( \sup_{\{t; e^{n-1} < v_t \leq e^n\}} \sup_{\{y,z \in K_\delta\}} |M_t^{y,z}| > \frac{u e^{n-1}}{M e^{n/4}}; \|M\|_{\rho, e^n} \leq \sqrt{M e^n} \right)$$

$$\leq 2 \exp\left( -\frac{u^2 e^{2n-2}}{M^3 e^{(3/2)n} \phi^2(\delta)} \right).$$

For all $u > 0$, for all $M > 0$ the series $\sum_n \mathbf{P}(B_n)$ converges, and the lemma follows by Borel–Cantelli. $\square$

Note that $v_t^{1/4}$ is enough for our needs, but actually one can replace $v_t^{1/4}$ in the lemma with $v_t^p$ without changing the proof.

**Lemma 2.3.** *For every $M > 0$ there is some $t = t(\omega, M)$ such that a.s. for $t > t(\omega, M)$ it holds*

$$\sup_{y \neq z \in K_\delta} \frac{1}{z - y} \int_y^z \frac{L_t^x}{V_t} \, dx \, \mathbf{1}_{\Omega_{M,t}} \leq A$$

*with some constant $A > 0$, independent of $M, t, \omega$.*

**Proof.** Let $y < z$ in $K_\delta$.

$$\int_y^z \frac{L_t^x}{V_t} \, dx \, \mathbf{1}_{\Omega_{M,t}} \leq (z - y) \left\{ \frac{L_t^{x_0}}{V_t} + \sup_{y,z} \frac{|L_t^y - L_t^z|}{V_t} \mathbf{1}_{\Omega_{M,t}} \right\}. \tag{12}$$

Using the Chacon–Ornstein theorem and Proposition 2.1, almost surely $\frac{L_t^{x_0}}{V_t} \to \frac{\|\nu_{L^x}\|}{\|\nu_V\|} = \sigma^2(x)\mu(x) := l_0$ with $l_0$ independent of $\omega$. Using the definition (7) and the notations (10) and (8) we can write:

$$\frac{|L_t^y - L_t^z|}{V_t} \mathbf{1}_{\Omega_{M,t}} \leq \frac{2\delta}{V_t} + \frac{|M_t^{y,z}|}{V_t} \mathbf{1}_{\Omega_{M,t}} + \frac{A_t}{V_t}.$$

Following Lemma 2.2 there is some $t = t(\omega, M)$, such that for $t > t(\omega, M)$

$$\frac{|L_t^y - L_t^z|}{V_t} \mathbf{1}_{\Omega_{M,t}} \leq \frac{2\delta}{V_t} + \frac{\delta}{v_t^{1/4}} + \frac{A_t}{V_t} \leq 2l_1, \tag{13}$$

where we have denoted $l_1 = \lim_{t \to \infty} \frac{A_t}{V_t}$. From (12) and (13) we deduce that the lemma holds with $A = 2(l_0 + l_1)$. $\square$

**Lemma 2.4.** *For every $M$ there is some $t = t(\omega, M)$, such that a.s. for $t > t(\omega, M)$ it holds with some constant $C > 0$*

$$\sup_{y \neq z \in K_\delta} \frac{|L_t^y - L_t^z|}{V_t} \mathbf{1}_{\Omega_{M,t}} \leq \delta C.$$

**Proof.** We apply Lemmas 2.2 and 2.3 to the decomposition

$$\frac{|L_t^y - L_t^z|}{V_t} \mathbf{1}_{\Omega_{M,t}} \leq \frac{2\delta}{V_t} + \frac{|M_t^{y,z}|}{V_t} \mathbf{1}_{\Omega_{M,t}} + \int_y^z \frac{L_t^x}{V_t} \, dx \, \mathbf{1}_{\Omega_{M,t}},$$

where in Lemma 2.2 we take $u = \delta$. $\square$

This property is weaker than the restricted uniform continuity of normalized local time, but enough for what we want to do. Namely, it provides the necessarily "uniform" argument to show that $\frac{\inf_K L_t^x}{v_t}$ is not far from $\frac{L_t^{x_0}}{v_t}$ for $K = [x_0 - \delta; x_0 + \delta]$ with a small $\delta$. Now we can prove Theorem 3:

**Proof.** We firstly prove the theorem for $K_\delta = [x_0 - \delta; x_0 + \delta]$ with some small $\delta$ depending on $x_0$.

Let $t = t(\omega, M)$ of Lemma 2.4. As previously, let $l_0 > 0$ be the almost sure limit $l_0 = \lim_{t \to \infty} \frac{L_t^{x_0}}{V_t}$. Let $t'(\omega) = \inf\{t > 0; l_0/2 < L_t^{x_0}/V_t < 2l_0\}$ and $\tau = \max(t, t')$. Fix $M > 0$ and $\delta = \delta(x_0) > 0$ in such a way that the following holds with $C > 0$ from Lemma 2.4

$$\frac{1}{M} \leq \frac{l_0}{2} - \delta C \leq 2l_0 + \delta C \leq M.$$

We have

$$\mathbf{P}\left(\frac{1}{M} \leq \inf_{x \in K_\delta} \frac{L_t^x}{V_t} \leq \sup_{x \in K_\delta} \frac{L_t^x}{V_t} \leq M\right)$$

$$\geq \mathbf{P}\left(\frac{1}{M} \leq \frac{L_t^{x_0}}{V_t} + \inf_{x \in K_\delta} \frac{L_t^x - L_t^{x_0}}{V_t} \leq \frac{L_t^{x_0}}{V_t} + \sup_{x \in K_\delta} \frac{L_t^x - L_t^{x_0}}{V_t} \leq M \cap \Omega_{M,t} \cap t > \tau\right)$$

$$\geq \mathbf{P}\left(\frac{1}{M} \leq \frac{L_t^{x_0}}{V_t} - \delta C \leq \frac{L_t^{x_0}}{V_t} + \delta C \leq M \cap \Omega_{M,t} \cap t > \tau\right)$$

$$\geq \mathbf{P}\left(\frac{1}{M} \leq \frac{l_0}{2} - \delta C \leq 2l_0 + \delta C \leq M \cap \Omega_{M,t} \cap t > \tau\right)$$

$$\geq \mathbf{P}(\Omega_{M,t} \cap t > \tau). \tag{14}$$

For the last expression we have

$$\lim_{M \to \infty} \liminf_{t \to \infty} \mathbf{P}(\Omega_{M,t} \cap t > \tau(\omega, M)) = 1.$$

Note that we can easily change $V_t$ with $v_t$ in (14), so the theorem for $K_\delta$ is proven.

Now let $K$ be some compact in the state space. For each $x \in K$ we choose $\delta_x > 0$ such that the theorem with $K(x) = [x - \delta_x; x + \delta_x]$ is true. Then by standard compactness arguments we obtain the theorem for $K$. □

The strong Chacon–Ornstein (SCO) theorem states that if $A_t$ is a SAF and $v_t$ is a deterministic equivalent associated to $X_t$, then $\mathbb{E}_\nu A_t \sim v_t$, as $t \to \infty$ for any initial distribution $\nu$. The following theorem is a uniform version of SCO theorem for a local time.

**Theorem 4.** *The sequence of real function* $g_t(y) = \frac{\mathbb{E}_\nu L_t^y}{v_t}$ *converges uniformly on $K_\delta$ to the function* $l(y) = \sigma^2(y)\mu(y)$ *as $t \to \infty$.*

**Proof.** The point by point convergence follows from the strong Chacon–Ornstein theorem and Proposition 2.1:

$$\lim_{t \to \infty} \frac{\mathbb{E}_\nu L_t^y}{v_t} = \frac{\|\nu_{L^y}\|}{1} = \sigma^2(y)\mu(y). \tag{15}$$

To show the uniform convergence we will prove that the family $\{\frac{\mathbb{E}_\nu L_t^y}{v_t}\}_{t>0}$ is equicontinuous on $K_\delta = [x_0 - \delta, x_0 + \delta]$.

$$\frac{1}{2} L_t^y = (X_t - y)^+ - (X_0 - y)^+ - \int_0^t \mathbf{1}_{\{X_s > y\}} \, dX_s$$

hence

$$\frac{1}{2}(\mathbb{E}_\nu L_t^y - \mathbb{E}_\nu L_t^z) = \mathbb{E}_\nu((X_t - y)^+ - (X_t - z)^+) - \mathbb{E}_\nu((X_0 - y)^+ - (X_0 - z)^+)$$

$$- \mathbb{E}_\nu \int_0^t \mathbf{1}_{\{y < X_s \leq z\}} b(X_s) \, ds.$$

So, as an initial inequality we can take

$$\frac{1}{2}|\mathbb{E}_\nu L_t^y - \mathbb{E}_\nu L_t^z| \leq 2|y - z| + \mathbb{E}_\nu \int_0^t \mathbf{1}_{\{y < X_s \leq z\}} |b(X_s)| \, ds. \tag{16}$$

We rewrite this using occupation formula, and divide by $v_t$:

$$\frac{|\mathbb{E}_\nu L_t^y - \mathbb{E}_\nu L_t^z|}{2v_t} \leq \frac{2|y-z|}{v_t} + \int_y^z \frac{|b(x)|}{\sigma^2(x)} \frac{\mathbb{E}_\nu L_t^x}{v_t} \, \mathrm{d}x. \tag{17}$$

From (16) we have for all $y, z \in K_\delta$

$$\frac{1}{2}|\mathbb{E}_\nu L_t^y - \mathbb{E}_\nu L_t^z| \leq 4\delta + \mathbb{E}_\nu \int_0^t \mathbf{1}_{\{X_s \in K_\delta\}} |b(X_s)| \, \mathrm{d}s = 4\delta + \tilde{v}_t. \tag{18}$$

Since $b$ is locally bounded and $X$ is strong Feller, the second term on the right-hand side in (18) is an expectation of special additive functional which was denoted by $\tilde{v}_t$. Hence, for all $x \in K_\delta$

$$\frac{\mathbb{E}_\nu L_t^x}{v_t} \leq \frac{|\mathbb{E}_\nu L_t^x - \mathbb{E}_\nu L_t^{x_0}|}{v_t} + \frac{\mathbb{E}_\nu L_t^{x_0}}{v_t} \leq \frac{4\delta + \tilde{v}_t + \tilde{\tilde{v}}_t}{v_t},$$

where we denoted $\mathbb{E}_\nu L_t^{x_0}$ by $\tilde{\tilde{v}}_t$ to stress the fact that $\mathbb{E}_\nu L_t^{x_0}$ is an expectation of another special additive functional. According to the strong Chacon–Ornstein theorem the right-hand side expression converges to some constant, hence for some constant $M > 0$ and $t > 0$ large enough

$$\sup_{x \in [x_0 - \delta, x_0 + \delta]} \frac{\mathbb{E}_\nu L_t^x}{v_t} \leq M. \tag{19}$$

Now we insert (19) in (17): for some $K > 0$ and a large enough $t$

$$\frac{|\mathbb{E}_\nu L_t^y - \mathbb{E}_\nu L_t^z|}{v_t} \leq \frac{|y-z|}{v_t} + M \sup_{[x_0-\delta, x_0+\delta]} \frac{|b(x)|}{\sigma^2(x)} |y-z| \leq K|y-z|$$

and uniform convergence follows from the Ascoli lemma. □

**Corolarry 2.1.** *Let $h_t, t \geq 0$, be some deterministic function with range in $[0, \delta]$ and such that $\lim_{t \to \infty} h_t = 0$. For all $y \in [-1; 1]$ $\lim_{t \to \infty} \frac{\mathbb{E}_\nu L_t^{x_0 + h_t y}}{v_t} = \mu(x_0)\sigma^2(x_0)$.*

**Proof.** We use the theorem about continuous convergence, see [13], p. 194: A sequence of continuous functions $(f_n)$ defined on some metric compact $X$ converges uniformly to $f$ if and only if for any sequence $(x_n)$ the condition $\lim_{n \to \infty} x_n = x$ implies $\lim_{n \to \infty} f_n(x_n) = f(x)$. □

**Corolarry 2.2.** *Let $\phi$ be measurable bounded and compactly supported in $[-1; 1]$. Let $h: \mathbb{R}^+ \to \mathbb{R}^+$ be such that $\lim_{t \to \infty} h_t = 0$, and $\psi: \mathbb{R} \to \mathbb{R}^+$ continuous in $x_0$. Then for any initial probability $\nu$ the following limit holds:*

$$\lim_{t \to \infty} \frac{\mathbb{E}_\nu \int_0^t \phi((X_s - x_0)/h_t)\psi(X_s) \, \mathrm{d}s}{h_t v_t} = \psi(x_0)\mu(x_0) \int_\mathbb{R} \phi(y) \, \mathrm{d}y.$$

**Proof.** Using the occupation formula

$$\frac{\mathbb{E}_\nu \int_0^t \phi((X_s - x_0)/h_t)\psi(X_s) \, \mathrm{d}s}{h_t v_t} = \frac{\mathbb{E}_\nu \int_\mathbb{R} \phi((x - x_0)/h_t)(\psi(x)/\sigma^2(x))L_t^x \, \mathrm{d}x}{h_t v_t}$$

$$= \int_\mathbb{R} \phi(y) \frac{\psi(x_0 + h_t y)}{\sigma^2(x_0 + h_t y)} f_t(y) \, \mathrm{d}y,$$

where we have denoted $f_t(y) = \frac{\mathbb{E}_\nu L_t^{x_0 + h_t y}}{v_t}$. Using Corollary 2.1,

$$\forall y \in [-1; 1] \quad \lim_{t \to \infty} f_t(y) = \mu(x_0)\sigma^2(x_0).$$

Also using (19) for some $M > 0$ and $t$ large enough we obtain $\sup_{[-1,1]} f_t \leq M$. The dominated convergence theorem concludes the proof. □

## 3. Limit theorems for AF and martingales composed with deterministic or random bandwidth

Let $\phi \colon \mathbb{R} \to \mathbb{R}^+$ be of class $\mathcal{C}^1$, with support in $[-1, 1]$ and $\int_{-1}^{1} \phi(x) \, dx = 1$. Let $X$ be a Harris recurrent diffusion subject to the assumptions of (1). Denote for $h \geq 0$

$$\frac{1}{h} A_t^h = \frac{1}{h} \int_0^t \phi\left(\frac{X_s - x_0}{h}\right) ds$$

(for $h = 0$ see (21)), and for $h > 0$

$$M_t^h = \int_0^t \phi\left(\frac{X_s - x_0}{h}\right) \sigma(X_s) \, dW_s. \tag{20}$$

Since $\phi$ is $\mathcal{C}^1$, there is a continuous modification of $\frac{1}{h} A_t^h$ and if $\sigma$ is bounded, of $M_t^h$, $h \in ]0, \delta]$. We adopt these modifications in what follows. As before, for $\delta > 0$, $K_\delta = [x_0 - \delta, x_0 + \delta]$ and $v_t$ denotes a deterministic equivalent of $X$.

In this section we obtain asymptotic results for $A_t^h$ and $M_t^h$ composed with deterministic or random bandwidth $h_t$.

**Theorem 5.** *For any initial probability $\nu$,*

$$\lim_{N \to \infty} \liminf_{t \to \infty} \mathbf{P}_\nu \left( \frac{1}{N} < \inf_{h \in [0;\delta]} \frac{A_t^h}{h v_t} \leq \sup_{h \in [0;\delta]} \frac{A_t^h}{h v_t} < N \right) = 1.$$

**Proof.** By occupation formula, there is a $P$ negligible set $\Omega_0$, outside of which we have

$$\frac{1}{h} A_t^h = \frac{1}{h} \int_0^t \phi\left(\frac{X_s - x_0}{h}\right) ds = \frac{1}{h} \int_\mathbb{R} \phi\left(\frac{x - x_0}{h}\right) \frac{1}{\sigma^2(x)} L_t^x \, dx$$

$$= \int_{-1}^{1} \frac{\phi(y)}{\sigma^2(x_0 + hy)} L_t^{x_0 + hy} \, dy \tag{21}$$

for all $t > 0$ and every $h \in [0, \delta]$.

With some $0 < C < \infty$, outside of $\Omega_0$, we have

$$\frac{1}{C} \inf_{x \in K_\delta} L_t^x \leq \inf_{h \in [0;\delta]} \frac{1}{h} A_t^h \leq \sup_{h \in [0;\delta]} \frac{1}{h} A_t^h \leq C \sup_{x \in K_\delta} L_t^x$$

for all $t > 0$. Fix some $N > 0$.

$$\mathbf{P}\left( \frac{1}{N} < \inf_{h \in [0;\delta]} \frac{A_t^h}{h v_t} \leq \sup_{h \in [0;\delta]} \frac{A_t^h}{h v_t} < N \right) \geq \mathbf{P}\left( \frac{C}{N} < \inf_{x \in K_\delta} \frac{L_t^x}{v_t} \leq \sup_{x \in K_\delta} \frac{L_t^x}{v_t} < \frac{N}{C} \right)$$

and the lemma follows from Theorem 3. □

**Corolarry 3.1.** *Let $H_t(\omega), t \geq 0$, be some measurable process with range in $[0, \delta]$. Denote by $A_t^{H_t}$ the composition of $A_t^h(\omega)$ with $H_t(\omega)$. Then for any initial probability $\nu$,*

$$\lim_{N \to \infty} \liminf_{t \to \infty} \mathbf{P}\left( \frac{1}{N} < \frac{A_t^{H_t}}{H_t v_t} < N \right) = 1.$$

**Lemma 3.1.** *Let $h: \mathbb{R}^+ \to ]0, \delta]$ be such that $\lim_{t \to \infty} h_t = 0$. Let $(M_t^h)$ be the martingale defined by (20) and $M_t^{h_t}$ is a composition of this martingale with $h_t$. Then for any initial probability $\nu$,*

$$\lim_{N \to \infty} \liminf_{t \to \infty} \mathbf{P}_\nu \left( -N < \frac{M_t^{h_t}}{\sqrt{h_t v_t}} < N \right) = 1.$$

**Proof.** $(M_s^{h_t}; s \in [0, t])$ is a martingale, so we have by the Markov inequality and quadratic variation's definition

$$\mathbf{P}_\nu \left( \left| \frac{M_t^{h_t}}{\sqrt{h_t v_t}} \right| > N \right) \leq \frac{\mathbb{E}_\nu [M_t^{h_t}]^2}{N^2 h_t v_t} = \frac{\mathbb{E}_\nu \int_0^t \phi^2((X_s - x_0)/h_t) \sigma^2(X_s) \, ds}{N^2 h_t v_t}.$$

Lemma 3.1 follows from Corollary 2.2 and the fact that $\sigma > 0$. □

**Lemma 3.2.** *Let $\nu$ be some probability on $\mathbb{R}$, $(H_t)_{t \geq 0}$ be adapted random processes with range in $]0, \delta]$ and $h_t$ some measurable deterministic function with range in $]0, \delta]$, such that*

$$\lim_{N \to \infty} \liminf_{t \to \infty} \mathbf{P}_\nu \left( \frac{1}{N} < \frac{H_t}{h_t} < N \right) = 1.$$

*Then we have*

$$\lim_{N \to \infty} \limsup_{t \to \infty} \mathbf{P}_\nu \left( \left| \frac{M_t^{H_t}}{\sqrt{h_t v_t}} - \frac{M_t^{h_t}}{\sqrt{h_t v_t}} \right| > N \right) = 0.$$

**Proof.** Under the assumption of the lemma we can differentiate with respect to $h$ the stochastic integral $M_t^h$. For $h \in ]0, \delta]$ its derivative is given by:

$$\frac{\partial}{\partial h} M_t^h = \int_0^t \phi' \left( \frac{X_s - x_0}{h} \right) \sigma(X_s) \, dW_s = \frac{1}{h} \int_0^t \psi \left( \frac{X_s - x_0}{h} \right) \sigma(X_s) \, dW_s = \frac{1}{h} N_t^h,$$

where we denote $\psi(\frac{X_s - x_0}{h}) = -\phi'(\frac{X_s - x_0}{h}) \frac{X_s - x_0}{h}$ and $N_t^h = \int_0^t \psi(\frac{X_s - x_0}{h}) \sigma(X_s) \, dW_s$. Remark that $N_t^h$ and $M_t^h$ have the same structure in the sense that $\psi(x)$ is also supported in $[-1, 1]$ and of class $\mathcal{C}^1(\mathbb{R})$. We will use the representation

$$M_t^{H_t} - M_t^{h_t} = \int_{h_t}^{H_t} \frac{1}{h} N_t^h \, dh.$$

We have

$$\mathbf{P} \left( \left| \frac{M_t^{H_t}}{\sqrt{h_t v_t}} - \frac{M_t^{h_t}}{\sqrt{h_t v_t}} \right| > N, \frac{1}{\sqrt{N}} \leq \frac{H_t}{h_t} \leq \sqrt{N} \right)$$

$$\leq \mathbf{P} \left( \frac{1}{\sqrt{h_t v_t}} \int_{h_t/\sqrt{N}}^{h_t \sqrt{N}} \frac{1}{h} |N_t^h| \, dh > N \right) \leq \frac{1}{N} \frac{1}{\sqrt{h_t v_t}} \mathbb{E}_\nu \int_{h_t/\sqrt{N}}^{h_t \sqrt{N}} \frac{1}{h} |N_t^h| \, dh$$

$$= \frac{1}{N} \frac{1}{\sqrt{h_t v_t}} \int_{h_t/\sqrt{N}}^{h_t \sqrt{N}} \frac{1}{\sqrt{h}} \frac{\mathbb{E}_\nu |N_t^h|}{\sqrt{h}} \, dh. \qquad (22)$$

Using the Cauchy–Schwarz inequality, quadratic-variation and occupation formulas

$$\frac{\mathbb{E}_\nu |N_t^h|}{\sqrt{h}} \leq \sqrt{\frac{\mathbb{E}_\nu |N_t^h|^2}{h}} = \sqrt{\frac{\mathbb{E}_\nu \int_0^t \psi^2((X_s - x_0)/h) \sigma^2(X_s) \, ds}{h}}$$

$$= \sqrt{\frac{\mathbb{E}_\nu \int_\mathbb{R} \psi^2((x-x_0)/h)\sigma^2(x) L_t^x \, \mathrm{d}x}{h}} = \sqrt{\int_\mathbb{R} \psi^2(y)\sigma^2(x_0+hy)\mathbb{E}_\nu L_t^{x_0+hy} \, \mathrm{d}y}.$$

According to Theorem 4 the sequence of real function $g_t(y) \to \frac{\mathbb{E}_\nu L_t^y}{v_t}$ converges uniformly on $K_\delta$ to the function $l(y) = \mu(y)\sigma^2(y)$ as $t \to \infty$. By assumptions $l$ is bounded on $K_\delta$, so there is some $A > 0$ and some $t_0$ such that for all $t > t_0$

$$\int_\mathbb{R} \psi^2(y)\sigma^2(x_0+hy)\frac{\mathbb{E}_\nu L_t^{x_0+hy}}{v_t} \, \mathrm{d}y \le A.$$

We insert now this estimation in (22):

$$\mathbf{P}_\nu\left(\left|\frac{M_t^{H_t}}{\sqrt{h_t v_t}} - \frac{M_t^{h_t}}{\sqrt{h_t v_t}}\right| > N, \frac{1}{\sqrt{N}} \le \frac{H_t}{h_t} \le \sqrt{N}\right)$$

$$\le \frac{1}{N}\frac{1}{\sqrt{h_t}}\int_{h_t/\sqrt{N}}^{h_t\sqrt{N}} \frac{1}{\sqrt{h}} \sqrt{\int_\mathbb{R} \psi^2(y)\sigma^2(x_0+hy)\frac{\mathbb{E}_\nu L_t^{x_0+hy}}{v_t} \, \mathrm{d}y} \, \mathrm{d}h$$

$$\le \sqrt{A}\frac{1}{N}\frac{1}{\sqrt{h_t}}\int_{h_t/\sqrt{N}}^{h_t\sqrt{N}} \frac{1}{\sqrt{h}} \, \mathrm{d}h = 2\sqrt{A}\frac{1}{N}\frac{1}{\sqrt{h_t}}\left(\sqrt{h_t}\sqrt{N} - \frac{\sqrt{h_t}}{\sqrt{N}}\right)$$

and the lemma follows by taking limit as $N \to \infty$.                                                                                                        □

From Lemmas 3.1 and 3.2 we deduce Corollary 3.1.

**Corolarry 3.2.** *Let $(H_t)_{t \ge 0}$ be adapted random processes with range in $]0, \delta]$. For any initial probability $\nu$,*

$$\lim_{N \to \infty} \liminf_{t \to \infty} \mathbf{P}_\nu\left(-N < \frac{M_t^{H_t}}{\sqrt{h_t v_t}} < N\right) = 1.$$

## 4. Application to non-parametric estimation

We observe until time $t$ a sample path of the processes $X$, weak solution of (1). Suppose that the drift coefficient $b$ is unknown and for $x_0 \in \mathbb{R}$ we want to estimate $b(x_0)$. Let $\delta > 0$ be such that $[x_0 - \delta; x_0 + \delta] \subset \mathbb{R}$. We suppose that $\sigma$ and $b$ satisfy the assumptions of (1), and also that for some $\alpha \in ]0, 1]$ and $\gamma > 0$ the function $b$ satisfies the local Hölder condition:

$$\sup_{x \in [x_0-\delta, x_0+\delta]} \frac{|b(x) - b(x_0)|}{|x - x_0|^\alpha} \le \gamma.$$

Let $\phi : \mathbb{R} \to \mathbb{R}^+$ be of class $C^1$ with support in $[-1, 1]$ and $\int \phi(x) \, \mathrm{d}x = 1$. For $h > 0$ we consider a family of Nadaraya–Watson estimators

$$\hat{b}_{x_0, t}^h = \frac{\int_0^t \phi((X_s - x_0)/h) \, \mathrm{d}X_s}{\int_0^t \phi((X_s - x_0)/h) \, \mathrm{d}s}.$$

Let $V_t$ be some given IAF of $X$ such that $\|\nu_V\| > 0$, and let $v_t$ be its deterministic equivalent.
Denote

$$H_t = V_t^{-1/(2\alpha+1)} \wedge \delta; \qquad h_t = v_t^{-1/(2\alpha+1)} \wedge \delta; \qquad r_t = v_t^{\alpha/(2\alpha+1)}; \qquad R_t = V_t^{\alpha/(2\alpha+1)}. \qquad (23)$$

To estimate $b(x_0)$ we use $\hat{b}_{x_0, t}^{H_t}$, so $H_t$ and $R_t$ play the role of "random bandwidth" and of "random rate," respectively.

**Theorem 6.** $R_t$ *is an upper rate of convergence of* $\hat{b}_{x_0,t}^{H_t}$ *to* $b(x_0)$

$$\lim_{K \to \infty} \limsup_{t \to \infty} \mathbf{P}_\nu(R_t|\hat{b}_{x_0,t}^{H_t} - b(x_0)| > K) = 0.$$

**Proof.** Using (1) we can write:

$$|\hat{b}_{x_0,t}^{H_t} - b(x_0)| \leq \frac{|\int_0^t \phi((X_s - x_0)/H_t)\sigma(X_s)\,\mathrm{d}W_s|}{\int_0^t \phi((X_s - x_0)/H_t)\,\mathrm{d}s} + \frac{\int_0^t \phi((X_s - x_0)/H_t)|b(X_s) - b(x_0)|\,\mathrm{d}s}{\int_0^t \phi((X_s - x_0)/H_t)\,\mathrm{d}s}$$

$$\leq \frac{|\int_0^t \phi((X_s - x_0)/H_t)\sigma(X_s)\,\mathrm{d}W_s|}{\int_0^t \phi((X_s - x_0)/H_t)\,\mathrm{d}s} + \gamma H_t^\alpha$$

$$= \frac{|\int_0^t \phi((X_s - x_0)/H_t)\sigma(X_s)\,\mathrm{d}W_s|}{\sqrt{v_t H_t}} \times \frac{v_t H_t}{\int_0^t \phi((X_s - x_0)/H_t)\,\mathrm{d}s} \times \frac{1}{\sqrt{v_t H_t}} + \gamma H_t^\alpha,$$

where in the second line we have used the Hölder property of $b$. Let

$$\Omega_{t,K}^1 = \left\{ \frac{|\int_0^t \phi((X_s - x_0)/H_t)\sigma(X_s)\,\mathrm{d}W_s|}{\sqrt{v_t H_t}} \leq K^{1/4} \right\};$$

$$\Omega_{t,K}^2 = \left\{ \frac{1}{K^{1/4}} \leq \frac{\int_0^t \phi((X_s - x_0)/H_t)\,\mathrm{d}s}{v_t H_t} \leq K^{1/4} \right\};$$

$$\Omega_{t,K}^3 = \left\{ \frac{1}{K^{1/4}} \leq \frac{V_t}{v_t} \leq K^{1/4} \right\}$$

and

$$\Omega_{t,K} = \Omega_{t,K}^1 \cap \Omega_{t,K}^2 \cap \Omega_{t,K}^3.$$

We have

$$\mathbf{P}_\nu(R_t \times |\hat{b}_{x_0,t}^{H_t} - b(x_0)| \geq K) \leq \mathbf{P}_\nu\left(v_t^{\alpha/(2\alpha+1)}\left(K^{1/2}\frac{1}{\sqrt{v_t h_t}} + \gamma h_t^\alpha\right) > K\right) + \mathbf{P}_\nu(\Omega_{t,K}^c)$$

$$= \mathbf{P}_\nu(K^{1/2} + 1 > K) + \mathbf{P}_\nu(\Omega_{t,K}^c),$$

where we used the fact that $h_t^\alpha = v_t^{-\alpha/(2\alpha+1)}$ for large $t$, $\sqrt{v_t h_t} = \sqrt{v_t^{1-1/(2\alpha+1)}} = v_t^{\alpha/(2\alpha+1)}$. Using Corollaries 3.1 and 3.2 we obtain

$$\lim_{K \to \infty} \limsup_{t \to \infty} \mathbf{P}(\Omega_{t,K}^c) = 0$$

and Theorem 6 follows. □

In the previous theorem we express the bandwidth and the rate as a function of additive functional $V_t$ of observed process $X_t$. It is easy to see that the theorem holds if we change $H_t$ with $h_t$ and $R_t$ with $r_t$ given by (23) which are both deterministic.

### Acknowledgments

We wish to thank S. Delattre and M. Hoffmann for useful comments and discussions, and the referee for his pertinent remarks.

## References

[1] A. Borodin and P. Salminen. *Handbook of Brownian Motion – Facts and Formulae*. Probability and its Applications. Birkhäuser, Basel, 1996. MR1477407

[2] M. Brancovan. Fonctionnelles additives spéciales des processus récurrents au sens de Harris. *Z. Wahrsch. Verw. Gebiete* **47** (1979) 163–194. MR0523168

[3] X. Chen. How often does a Harris recurrent Markov chain recur? *Ann. Probab.* **27** (1999) 1324–1346. MR1733150

[4] A. Dalalyan. Sharp adaptive estimation of the drift function for ergodic diffusions. *Ann. Statist.* **33** (2005) 2507–2528. MR2253093

[5] A. Dalalyan and Y. Kutoyants. On second order minimax estimation of the invariant density for ergodic diffusions. *Statist. Decisions* **22** (2004) 17–41. MR2065989

[6] S. Delattre and M. Hoffmann. Asymptotic equivalence for a null recurrent diffusion model. *Bernoulli* **8** (2002) 139–174. MR1895888

[7] S. Delattre, M. Hoffmann and M. Kessler. Dynamics adaptive estimation of a scalar diffusion. Prépublication PMA-762, Univ. Paris 6. Available at www.proba.jussieu.fr/mathdoc/preprints/.

[8] L. Galtchouk and S. Pergamentchikov. Sequential nonparametric adaptive estimation of the drift coefficient in diffusion processes. *Math. Methods Statist.* **10** (2001) 316–330. MR1867163

[9] R. Höpfner and Y. Kutoyants. On a problem of statistical inference in null recurrent diffusions. *Stat. Inference Stoch. Process.* **6** (2003) 25–42. MR1965183

[10] R. Höpfner and E. Löcherbach. *Limit Theorems for Null Recurrent Markov Processes*. Providence, RI, 2003. MR1949295

[11] K. Itô and H. P. McKean, Jr. *Diffusion Processes and Their Sample Paths.* Springer, Berlin, 1974. MR0345224

[12] R. Khasminskii. Limit distributions of some integral functionals for null-recurrent diffusions. *Stochastic Process. Appl.* **92** (2001) 1–9. MR1815176

[13] K. Kuratowski. *Introduction a la theorie des ensembles et a la topologie*. Institut de Mathematiques, Universite Geneve, 1966. MR0231338

[14] Y. Kutoyants. *Statistical Inference for Ergodic Diffusion Processes*. Springer, London, 2004. MR2144185

[15] E. Löcherbach and D. Loukianova. On Nummelin splitting for continuous time Harris recurrent Markov processes and application to kernel estimation for multidimensional diffusions. To appear in *Stochastic Process. Appl.*

[16] D. Loukianova and O. Loukianov. Deterministic equivalents of additive functionals of recurrent diffusions and drift estimation. To appear in *Stat. Inference Stoch. Process.*

[17] D. Loukianova and O. Loukianov. Almost sure rate of convergence of maximum likelihood estimators for multidimensional diffusions. In *Dependence in Probability and Statistics* 269–347. Springer, New York, 2006. MR2283262

[18] Y. Nishiyama. A maximum inequality for continuous martingales and $M$-estimation in Gaussian white noise model. *Ann. Statist.* **27** (1999) 675–696. MR1714712

[19] D. Revuz and M. Yor. *Continuous Martingales and Brownian Motion*. Springer, Berlin, 1994. MR1303781

[20] L. C. G. Rogers and D. Williams. *Diffusions, Markov Processes, and Martingales*, Vol. 2, Wiley, New York, 1990. MR0921238

[21] A. Touati. Théorèmes limites pour les processus de Markov récurrents. Unpublished paper, 1988. (See also *C.R.A.S. Paris Série I* **305** (1987) 841–844.) MR0923211

[22] H. van Zanten. On empirical processes for ergodic diffusions and rates of convergence of $M$-estimators. *Scand. J. Statist.* **30** (2003) 443–458. MR2002221

[23] H. van Zanten. On the rate of convergence of the maximum likelihood estimator in Brownian semimartingale models. *Bernoulli* **11** (2005) 643–664. MR2158254

[24] N. Yoshida. Asymptotic behavior of $M$-estimators and related random field for diffusion process. *Ann. Inst. Statist. Math.* **42** (1990) 221–251. MR1064786